\documentclass{amsproc}%
\usepackage{amsmath}
\usepackage{amsfonts}
\usepackage{amssymb}
\usepackage{graphicx}%
\providecommand{\U}[1]{\protect\rule{.1in}{.1in}}
\newtheorem{theorem}{Theorem}[section]
\newtheorem{lemma}[theorem]{Lemma}

\newtheorem{corollary}[theorem]{Corollary}

\newtheorem{proposition}[theorem]{Proposition}

\theoremstyle{definition}

\theoremstyle{remark}

\numberwithin{equation}{section}

\begin{document}
\title{Some properties of the inverse error function}
\author{Diego Dominici }
\address{Department of Mathematics\\
State University of New York at New Paltz\\
75 S. Manheim Blvd. Suite 9\\
New Paltz, NY 12561-2443\\
USA\\
Phone: (845) 257-2607\\
Fax: (845) 257-3571}
\email{dominicd@newpaltz.edu}
\thanks{This work was partially supported by a Provost Research Award from SUNY New Paltz.}
\subjclass{Primary 33B20; Secondary 30B10, 34K25}
\date{June 4, 2007}
\keywords{Inverse error function, asymptotic analysis, discrete ray method,
differential-difference equations, Taylor series}

\begin{abstract}
The inverse of the error function, $\operatorname{inverf}(x),$ has
applications in diffusion problems, chemical potentials, ultrasound imaging,
etc. We analyze the derivatives $\left.  \frac{d^{n}}{dz^{n}}%
\operatorname*{inverf}\left(  z\right)  \right\vert _{z=0}$, as $n\rightarrow
\infty$ using nested derivatives and a discrete ray method. We obtain a very
good approximation of $\operatorname{inverf}(x)$ through a high-order Taylor
expansion around $x=0$. We give numerical results showing the accuracy of our formulas.

\end{abstract}
\maketitle






\section{Introduction}

The error function $\operatorname{erf}(z),$ defined by
\[
\operatorname{erf}(z)=\frac{2}{\sqrt{\pi}}\int\limits_{0}^{z}\exp\left(
-t^{2}\right)  dt,
\]
occurs widely in almost every branch of applied mathematics and mathematical
physics, e.g., probability and statistics \cite{MR0034250}, data analysis
\cite{MR999553}, heat conduction \cite{MR0016873}, etc. It plays a fundamental
role in asymptotic expansions \cite{MR1429619} and exponential asymptotics
\cite{MR990851}.

Its inverse, which we will denote by $\operatorname*{inverf}\left(  z\right)
,$
\[
\operatorname*{inverf}\left(  z\right)  =\operatorname{erf}^{-1}(z),
\]
appears in multiple areas of mathematics and the natural sciences. A few
examples include concentration-dependent diffusion problems \cite{MR0071876},
\cite{MR0281322}, solutions to Einstein's scalar-field equations
\cite{PhysRevD.51.444}, chemical potentials \cite{MR2166352}, the distribution
of lifetimes in coherent-noise models \cite{PhysRevE.59.R2512}, diffusion
rates in tree-ring chemistry \cite{MR2142222} and $3D$ freehand ultrasound
imaging \cite{san-joseMICCAI03}.

Although some authors have studied the function $\operatorname*{inverf}\left(
z\right)  $ (see \cite{MR1986919} and references therein), little is known
about its analytic properties$,$ the major work having been done in developing
algorithms for numerical calculations \cite{MR0341812}. Dan Lozier, remarked
the need for new techniques in the computation of $\operatorname*{inverf}%
\left(  z\right)  $ \cite{MR1393742}.

In this paper, we analyze the asymptotic behavior of the derivatives $\left.
\frac{d^{n}}{dz^{n}}\operatorname*{inverf}\left(  z\right)  \right\vert
_{z=0}$ for large values of $n,$ using a discrete WKB method \cite{MR1373150}.
In Section 2 we present some properties of the derivatives of
$\operatorname*{inverf}\left(  z\right)  $ and review our previous work on
nested derivatives. In Section 3 we study a family of polynomials $P_{n}(x)$
associated with the derivatives of $\operatorname*{inverf}\left(  z\right)  $,
which were introduced by L. Carlitz in \cite{MR0153878}. Theorem \ref{theorem}
contains our main result on the asymptotic analysis of $P_{n}(x).$ In Section
4 we give asymptotic approximations for $\left.  \frac{d^{n}}{dz^{n}%
}\operatorname*{inverf}\left(  z\right)  \right\vert _{z=0}$ and some
numerical results testing the accuracy of our formulas.

\section{Derivatives}

Let us denote the function $\operatorname*{inverf}\left(  z\right)  $ by
$\mathfrak{I}(z)$ and its derivatives by
\begin{equation}
d_{n}=\left.  \frac{d^{n}}{dz^{n}}\operatorname*{inverf}\left(  z\right)
\right\vert _{z=0},\quad n=0,1,\ldots. \label{dn}%
\end{equation}
Since $\operatorname{erf}(z)$ tends to $\pm1$ as $z\rightarrow\pm\infty,$ it
is clear that $\operatorname*{inverf}\left(  z\right)  $ is defined in the
interval $\left(  -1,1\right)  $ and has singularities at the end points.

\begin{proposition}
The function $\mathfrak{I}(z)$ satisfies the nonlinear differential equation%
\begin{equation}
\mathfrak{I}^{\prime\prime}-2\mathfrak{I}\left(  \mathfrak{I}^{\prime}\right)
^{2}=0 \label{ODE}%
\end{equation}
with initial conditions%
\begin{equation}
\mathfrak{I}(0)=0,\quad\mathfrak{I}^{\prime}(0)=\frac{\sqrt{\pi}}{2}.
\label{d0d1}%
\end{equation}

\end{proposition}

\begin{proof}
It is clear that $\mathfrak{I}(0)=0,$ since $\operatorname{erf}(0)=0.$ Using
the chain rule, we have
\[
\mathfrak{I}^{\prime}\left[  \operatorname{erf}(z)\right]  =\frac
{1}{\operatorname{erf}^{\prime}(z)}=\frac{\sqrt{\pi}}{2}\exp\left\{
\mathfrak{I}^{2}\left[  \operatorname{erf}(z)\right]  \right\}
\]
and therefore%
\begin{equation}
\mathfrak{I}^{\prime}=\frac{\sqrt{\pi}}{2}\exp\left(  \mathfrak{I}^{2}\right)
. \label{I'}%
\end{equation}
Setting $z=0$ we get $\mathfrak{I}^{\prime}(0)=\frac{\sqrt{\pi}}{2}$ and
taking the logarithmic derivative of (\ref{I'}) the result follows.
\end{proof}

To compute higher derivatives of $\mathfrak{I}(z),$ we begin by establishing
the following corollary.

\begin{corollary}
The function $\mathfrak{I}(z)$ satisfies the nonlinear differential-integral
equation%
\begin{equation}
\mathfrak{I}^{\prime}(z)\int\limits_{0}^{z}\mathfrak{I}(t)dt=-\frac{1}%
{2}+\frac{1}{\sqrt{\pi}}\mathfrak{I}^{\prime}(z). \label{int diff}%
\end{equation}

\end{corollary}

\begin{proof}
Rewriting (\ref{ODE}) as%
\[
\mathfrak{I}=\frac{1}{2}\frac{\mathfrak{I}^{\prime\prime}}{\left(
\mathfrak{I}^{\prime}\right)  ^{2}}%
\]
and integrating, we get%
\[
\int\limits_{0}^{z}\mathfrak{I}(t)dt=\frac{1}{2}\left[  -\frac{1}%
{\mathfrak{I}^{\prime}(z)}+\frac{1}{\mathfrak{I}^{\prime}(0)}\right]
=\frac{1}{2}\left[  -\frac{1}{\mathfrak{I}^{\prime}(z)}+\frac{2}{\sqrt{\pi}%
}\right]
\]
and multiplying by $\mathfrak{I}^{\prime}(z)$ we obtain (\ref{int diff}).
\end{proof}

\begin{proposition}
The derivatives of $\mathfrak{I}(z)$ satisfy the nonlinear recurrence%
\begin{equation}
d_{n+1}=\sqrt{\pi}\sum\limits_{k=0}^{n-1}\binom{n}{k+1}d_{k}d_{n-k},\quad
n=1,2,\ldots\label{recurrence}%
\end{equation}
with $d_{0}=0$ and $d_{1}=\frac{\sqrt{\pi}}{2}.$
\end{proposition}

\begin{proof}
Using
\[
\mathfrak{I}(z)=\sum\limits_{n=0}^{\infty}d_{n}\frac{z^{n}}{n!}%
\]
and $d_{1}=\frac{\sqrt{\pi}}{2}$ in (\ref{int diff}), we have%
\[
\left[  \frac{\sqrt{\pi}}{2}+\sum\limits_{n=1}^{\infty}d_{n+1}\frac{z^{n}}%
{n!}\right]  \left[  \sum\limits_{n=1}^{\infty}d_{n-1}\frac{z^{n}}{n!}%
-\frac{1}{\sqrt{\pi}}\right]  =-\frac{1}{2}%
\]
or%
\[
\frac{\sqrt{\pi}}{2}\sum\limits_{n=1}^{\infty}d_{n-1}\frac{z^{n}}{n!}%
+\sum\limits_{n=2}^{\infty}\left[  \sum\limits_{k=0}^{n-2}\binom{n}{k+1}%
d_{k}d_{n-k}\right]  \frac{z^{n}}{n!}-\frac{1}{\sqrt{\pi}}\sum\limits_{n=1}%
^{\infty}d_{n+1}\frac{z^{n}}{n!}=0.
\]
Comparing powers of $z^{n},$ we get%
\[
\frac{\sqrt{\pi}}{2}d_{n-1}+\sum\limits_{k=0}^{n-2}\binom{n}{k+1}d_{k}%
d_{n-k}-\frac{1}{\sqrt{\pi}}d_{n+1}=0
\]
or%
\[
\sum\limits_{k=0}^{n-1}\binom{n}{k+1}d_{k}d_{n-k}-\frac{1}{\sqrt{\pi}}%
d_{n+1}=0.
\]

\end{proof}

Although one could use (\ref{recurrence}) to compute the higher derivatives of
$\operatorname*{inverf}\left(  z\right)  ,$ the nonlinearity of the recurrence
makes it hard to analyze the asymptotic behavior of $d_{n}$ as $n\rightarrow
\infty.$ Instead, we shall use an alternative technique that we developed in
\cite{MR2031140} and we called the method of "nested derivatives". The
following theorem contains the main result presented in \cite{MR2031140}.

\begin{theorem}
Let
\[
H(x)=h^{-1}(x),\quad f(x)=\frac{1}{h^{\prime}(x)},\quad z_{0}=h(x_{0}),\text{
\ }\ \left\vert f(x_{0})\right\vert \in\left(  0,\infty\right)  .
\]
$\ $ \ Then,%
\[
H(z)=x_{0}+f(x_{0})\sum\limits_{n=1}^{\infty}\mathfrak{D}^{n-1}[f]\,(x_{0}%
)\frac{(z-z_{0})^{n}}{n!},
\]
where we define $\mathfrak{D}^{n}[f]$\thinspace$(x),$ \textit{the n}$^{th}%
$\textit{ nested derivative} of the function $f(x),$ by $\mathfrak{D}%
^{0}[f]\,(x)=1$ and%
\begin{equation}
\mathfrak{D}^{n+1}[f]\,(x)=\frac{d}{dx}\left[  f(x)\times\mathfrak{D}%
^{n}[f]\,(x)\right]  ,\quad n=0,1,\ldots. \label{nested}%
\end{equation}

\end{theorem}

The following proposition makes the computation of $\mathfrak{D}%
^{n-1}[f]\,(x_{0})$ easier in some cases.

\begin{proposition}
Let%
\begin{equation}
\mathfrak{D}^{n}[f]\,(x)=\sum\limits_{k=0}^{\infty}A_{k}^{n}\frac
{(x-x_{0})^{k}}{k!},\qquad f(x)=\sum\limits_{k=0}^{\infty}B_{k}\frac
{(x-x_{0})^{k}}{k!}. \label{AB}%
\end{equation}
Then,
\begin{equation}
A_{k}^{n+1}=\left(  k+1\right)  \sum\limits_{j=0}^{k+1}A_{k+1-j}^{n}B_{j}.
\label{A}%
\end{equation}

\end{proposition}

\begin{proof}
From (\ref{AB}) we have%
\begin{equation}
f(x)\mathfrak{D}^{n}[f]\,(x)=\sum\limits_{k=0}^{\infty}\alpha_{k}^{n}%
\frac{(x-x_{0})^{k}}{k!}, \label{Dfxf}%
\end{equation}
with%
\begin{equation}
\alpha_{k}^{n}=\sum\limits_{j=0}^{k}A_{k-j}^{n}B_{j}. \label{alpha}%
\end{equation}
Using (\ref{AB}) and (\ref{Dfxf}) in (\ref{nested}), we obtain%
\[
\sum\limits_{k=0}^{\infty}A_{k}^{n+1}(x-x_{0})^{k}=\frac{d}{dx}\sum
\limits_{k=0}^{\infty}\alpha_{k}^{n}(x-x_{0})^{k}=\sum\limits_{k=0}^{\infty
}\left(  k+1\right)  \alpha_{k+1}^{n}(x-x_{0})^{k}%
\]
and the result follows from (\ref{alpha}).
\end{proof}

To obtain a linear relation between successive nested derivatives, we start by
establishing the following lemma.

\begin{lemma}
Let
\begin{equation}
g_{n}\left(  x\right)  =\frac{\mathfrak{D}^{n}[f]\,(x)}{f^{n}\left(  x\right)
}. \label{gn}%
\end{equation}
Then,%
\begin{equation}
g_{n+1}(x)=g_{n}^{\prime}\left(  x\right)  +\left(  n+1\right)  \frac
{f^{\prime}(x)}{f(x)}g_{n}\left(  x\right)  ,\quad n=0,1,\ldots.
\label{ddnested}%
\end{equation}

\end{lemma}

\begin{proof}
Using (\ref{nested}) in (\ref{gn}), we have%
\begin{gather*}
g_{n+1}\left(  x\right)  =\frac{\mathfrak{D}^{n+1}[f]\,(x)}{f^{n+1}\left(
x\right)  }=\frac{\frac{d}{dx}\left[  f(x)\times\mathfrak{D}^{n}%
[f]\,(x)\right]  }{f^{n+1}\left(  x\right)  }\\
=\frac{\frac{d}{dx}\left[  g_{n}\left(  x\right)  f^{n+1}\left(  x\right)
\right]  }{f^{n+1}\left(  x\right)  }=\frac{g_{n}^{\prime}\left(  x\right)
f^{n+1}\left(  x\right)  +g_{n}\left(  x\right)  (n+1)f^{n}\left(  x\right)
f^{\prime}(x)}{f^{n+1}\left(  x\right)  }%
\end{gather*}
and the result follows.
\end{proof}

\begin{corollary}
Let
\[
H(x)=h^{-1}(x),\quad f(x)=\frac{1}{h^{\prime}(x)},\quad z_{0}=h(x_{0}),\text{
\ }\ \left\vert f(x_{0})\right\vert \in\left(  0,\infty\right)  .
\]
Then,$\ $
\begin{equation}
\frac{d^{n}H}{dz^{n}}(z_{0})=\left[  f(x_{0})\right]  ^{n}g_{n-1}(x_{0}),\quad
n=1,2,\ldots. \label{deriv}%
\end{equation}

\end{corollary}

For the function $h(x)=\operatorname{erf}(z),$ we have
\begin{equation}
f(x)=\frac{1}{h^{\prime}(x)}=\frac{\sqrt{\pi}}{2}\exp\left(  x^{2}\right)  ,
\label{f}%
\end{equation}
and setting $x_{0}=0$ we obtain $z_{0}=\operatorname{erf}(0)=0.$ Using the
Taylor series%
\[
\frac{\sqrt{\pi}}{2}\exp\left(  x^{2}\right)  =\frac{\sqrt{\pi}}{2}%
\sum\limits_{k=0}^{\infty}\frac{x^{2k}}{k!}%
\]
in (\ref{A}), we get%
\[
A_{k}^{n+1}=\frac{\sqrt{\pi}}{2}\left(  k+1\right)  \sum\limits_{j=0}%
^{\left\lfloor \frac{k+1}{2}\right\rfloor }\frac{A_{k+1-2j}^{n}}{j!},
\]
with $A_{k}^{n}$ defined in (\ref{AB}). Using (\ref{f}) in (\ref{ddnested}),
we have%
\begin{equation}
g_{n+1}(x)=g_{n}^{\prime}\left(  x\right)  +2\left(  n+1\right)  xg_{n}\left(
x\right)  ,\quad n=0,1,\ldots, \label{g}%
\end{equation}
while (\ref{deriv}) gives%
\begin{equation}
d_{n}=\left(  \frac{\sqrt{\pi}}{2}\right)  ^{n}g_{n-1}(0),\quad n=1,2,\ldots.
\label{deriv2}%
\end{equation}

In the next section we shall find an asymptotic approximation for a family of
polynomials closely related to $g_{n}\left(  x\right)  $.

\section{The polynomials $P_{n}(x)$}

We define the polynomials $P_{n}(x)$ by $P_{0}(x)=1$ and
\begin{equation}
P_{n}(x)=g_{n}\left(  \frac{x}{\sqrt{2}}\right)  2^{-\frac{n}{2}}.
\label{Pngn}%
\end{equation}%
\begin{equation}
P_{n+1}(x)=P_{n}^{\prime}(x)+\left(  n+1\right)  xP_{n}(x), \label{diffdiff}%
\end{equation}
The first few $P_{n}\left(  x\right)  $ are
\[
P_{1}(x)=x,\quad P_{2}(x)=1+2x^{2},\quad P_{3}(x)=7x+6x^{3},~\ldots~.
\]

The following propositions describe some properties of $P_{n}\left(  x\right)
.$

\begin{proposition}
Let%
\begin{equation}
P_{n}(x)=\sum\limits_{k=0}^{\left\lfloor \frac{n}{2}\right\rfloor }C_{k}%
^{n}x^{n-2k}, \label{Pcnk}%
\end{equation}
where $\left\lfloor \cdot\right\rfloor $ denotes the integer part function.
Then,%
\begin{equation}
C_{0}^{n}=n! \label{C01}%
\end{equation}
and%
\begin{equation}
C_{k}^{n}=n!\sum\limits_{j_{k}=0}^{n-1}\sum\limits_{j_{k-1}=0}^{j_{k}-1}%
\cdots\sum\limits_{j_{1}=0}^{j_{2}-1}\prod_{i=1}^{k}\frac{j_{i}-2i+2}{j_{i}%
+1},\quad k=1,\ldots,\left\lfloor \frac{n}{2}\right\rfloor . \label{cnk}%
\end{equation}

\end{proposition}

\begin{proof}
Using (\ref{Pcnk}) in (\ref{diffdiff}) we have%
\begin{gather*}
\sum\limits_{0\leq2k\leq n+1}^{{}}C_{k}^{n+1}x^{n+1-2k}=\sum\limits_{0\leq
2k\leq n}^{{}}C_{k}^{n}\left(  n-2k\right)  x^{n-2k-1}+\sum\limits_{0\leq
2k\leq n}^{{}}\left(  n+1\right)  C_{k}^{n}x^{n+1-2k}\\
=\sum\limits_{2\leq2k\leq n+2}^{{}}C_{k-1}^{n}\left(  n-2k+2\right)
x^{n+1-2k}+\sum\limits_{0\leq2k\leq n}^{{}}\left(  n+1\right)  C_{k}%
^{n}x^{n+1-2k}.
\end{gather*}
Comparing coefficients in the equation above, we get%
\begin{equation}
C_{0}^{n+1}=C_{0}^{n}, \label{C0}%
\end{equation}%
\begin{equation}
C_{k}^{n+1}=\left(  n-2k+2\right)  C_{k-1}^{n}+\left(  n+1\right)  C_{k}%
^{n},\quad k=1,\ldots,\left\lfloor \frac{n}{2}\right\rfloor \label{Cn+1}%
\end{equation}
and for $n=2m-1,$%
\[
C_{m}^{2m}=C_{m-1}^{2m-1},\quad m=1,2,\ldots.
\]
From (\ref{C0}) we immediately conclude that $C_{0}^{n}=n!,$ while
(\ref{Cn+1}) gives%
\begin{equation}
C_{k}^{n}=n!\sum\limits_{j=0}^{n-1}\frac{j-2k+2}{\left(  j+1\right)  !}%
C_{k-1}^{j},\quad n,k\geq1. \label{Cnk1}%
\end{equation}

Setting $k=1$ in (\ref{Cnk1}) and using (\ref{C01}), we have%
\begin{equation}
C_{1}^{n}=n!\sum\limits_{j=0}^{n-1}\frac{j}{\left(  j+1\right)  !}C_{0}%
^{j}=n!\sum\limits_{j=0}^{n-1}\frac{j}{j+1}. \label{C1n}%
\end{equation}
Similarly, setting $k=2$ in (\ref{Cnk1}) and using (\ref{C1n}), we get%
\[
C_{2}^{n}=n!\sum\limits_{j=0}^{n-1}\frac{j-2}{\left(  j+1\right)  !}\left[
j!\sum\limits_{i=0}^{j-1}\frac{i}{i+1}\right]  =n!\sum\limits_{j=0}^{n-1}%
\sum\limits_{i=0}^{j-1}\frac{j-2}{j+1}\frac{i}{i+1}%
\]
and continuing this way we obtain (\ref{cnk}).
\end{proof}

\begin{proposition}
The zeros of the polynomials $P_{n}(x)$ are purely imaginary for $n\geq1.$
\end{proposition}

\begin{proof}
For $n=1$ the result is obviously true. Assuming that it is true for $n$ and
that $P_{n}(x)$ is written in the form%
\begin{equation}
P_{n}(x)=n!%
{\displaystyle\prod\limits_{k=1}^{n}}
(z-z_{k}),\quad\operatorname{Re}(z_{k})=0,\quad1\leq k\leq n, \label{product}%
\end{equation}
we have two possibilities for $z^{\ast},$ with$\ P_{n+1}(z^{\ast})=0$:

\begin{enumerate}
\item $z^{\ast}=z_{k}$, for some $1\leq k\leq n.$

In this case, $\operatorname{Re}(z^{\ast})=0$ and the proposition is proved.

\item $z^{\ast}\neq z_{k}$, for all $1\leq k\leq n$.

From (\ref{diffdiff}) and (\ref{product}) we get%
\[
\frac{P_{n+1}(x)}{P_{n}(x)}=\frac{d}{dx}\ln\left[  P_{n}(x)\right]  +(n+1)x=%
{\displaystyle\sum\limits_{k=1}^{n}}
\frac{1}{z-z_{k}}+(n+1)x.
\]
Evaluating at $z=z^{\ast},$ we obtain%
\[
0=%
{\displaystyle\sum\limits_{k=1}^{n}}
\frac{1}{z^{\ast}-z_{k}}+(n+1)z^{\ast}%
\]
and taking $\operatorname{Re}(\bullet),$ we have%
\begin{gather*}
0=\operatorname{Re}\left[
{\displaystyle\sum\limits_{k=1}^{n}}
\frac{1}{z^{\ast}-z_{k}}+(n+1)z^{\ast}\right] \\
=%
{\displaystyle\sum\limits_{k=1}^{n}}
\frac{\operatorname{Re}\left(  z^{\ast}-z_{k}\right)  }{\left\vert z^{\ast
}-z_{k}\right\vert ^{2}}+(n+1)\operatorname{Re}(z^{\ast})=\operatorname{Re}%
(z^{\ast})\left[
{\displaystyle\sum\limits_{k=1}^{n}}
\frac{1}{\left\vert z^{\ast}-z_{k}\right\vert ^{2}}+n+1\right]
\end{gather*}
which implies that \ $\operatorname{Re}(z^{\ast})=0.$
\end{enumerate}
\end{proof}

\subsection{Asymptotic analysis of $P_{n}(x)$}

We first consider solutions to (\ref{diffdiff}) of the form%
\begin{equation}
P_{n}(x)=n!A^{\left(  n+1\right)  }(x), \label{Pn1}%
\end{equation}
with $x>0.$ Replacing (\ref{Pn1}) in (\ref{diffdiff}) and simplifying the
resulting expression, we obtain%
\[
A^{2}(x)=A^{\prime}(x)+xA(x),
\]
with solution%
\begin{equation}
A(x)=\exp\left(  -\frac{x^{2}}{2}\right)  \left[  C-\sqrt{\frac{\pi}{2}%
}\operatorname{erf}\left(  \frac{x}{\sqrt{2}}\right)  \right]  ^{-1},
\label{f1}%
\end{equation}
for some constant $C.$ Note that (\ref{Pn1}) is not an exact solution of
(\ref{diffdiff}), since it does not satisfy the initial condition $P_{0}(x)=1.$
To determine $C$ in(\ref{f1}), we observe from (\ref{C01}) that
\begin{equation}
P_{n}(x)\sim n!x^{n},\quad x\rightarrow\infty. \label{Pn4}%
\end{equation}
As $x\rightarrow\infty,$ we get from (\ref{f1})%
\[
\ln\left[  A(x)\right]  \sim-\frac{x^{2}}{2}-\ln\left(  C-\sqrt{\frac{\pi}{2}%
}\right)  +\frac{\exp\left(  -\frac{x^{2}}{2}\right)  }{\left(  C-\sqrt
{\frac{\pi}{2}}\right)  x},\quad x\rightarrow\infty,
\]
which is inconsistent with (\ref{Pn4}) unless $C=\sqrt{\frac{\pi}{2}}.$ In
this case, we have%
\begin{equation}
A(x)\sim x+\frac{1}{x},\quad x\rightarrow\infty, \label{A1}%
\end{equation}
matching (\ref{Pn4}). Thus,
\begin{equation}
A(x)=\sqrt{\frac{2}{\pi}}\exp\left(  -\frac{x^{2}}{2}\right)  \left[
1-\operatorname{erf}\left(  \frac{x}{\sqrt{2}}\right)  \right]  ^{-1}.
\label{psi2}%
\end{equation}
Since (\ref{Pn1}) and (\ref{A1}) give
\[
P_{n}(x)\sim n!x^{n+1},\quad x\rightarrow\infty,
\]
instead of (\ref{Pn4}), we need to consider
\begin{equation}
P_{n}(x)=n!A^{\left(  n+1\right)  }(x)B(x,n). \label{Pn2}%
\end{equation}
Replacing (\ref{Pn2}) in (\ref{diffdiff}) and simplifying, we get%
\[
B(x,n+1)=B(x,n)+\frac{1}{A(x)(n+1)}\frac{\partial B}{\partial x}(x,n).
\]
Using the approximation
\[
B(x,n+1)=B(x,n)+\frac{\partial B}{\partial n}(x,n)+\frac{1}{2}\frac
{\partial^{2}B}{\partial n^{2}}(x,n)+\cdots,
\]
we obtain%
\[
\frac{\partial B}{\partial n}=\frac{1}{A(x)(n+1)}\frac{\partial B}{\partial
x},
\]
whose solution is%
\begin{equation}
B(x,n)=F\left[  \frac{n+1}{1-\operatorname{erf}\left(  \frac{x}{\sqrt{2}%
}\right)  }\right]  , \label{B1}%
\end{equation}
for some function $F(u).$ Matching (\ref{Pn2}) with (\ref{Pn4}) requires
\begin{equation}
B(x,n)\sim\frac{1}{x},\quad x\rightarrow\infty. \label{B2}%
\end{equation}
Since in the limit as $x\rightarrow\infty,$ with $n$ fixed we have%
\[
\ln\left[  \frac{n+1}{1-\operatorname{erf}\left(  \frac{x}{\sqrt{2}}\right)
}\right]  \sim\frac{x^{2}}{2},
\]
(\ref{B1})-(\ref{B2}) imply%
\[
F(u)=\frac{1}{\sqrt{2\ln(u)}}.
\]
Therefore, for $x>0,$
\begin{equation}
P_{n}(x)\sim n!\Phi\left(  x,n\right)  ,\quad n\rightarrow\infty, \label{Pn5}%
\end{equation}
with%
\[
\Phi\left(  x,n\right)  =\left[  \sqrt{\frac{2}{\pi}}\frac{\exp\left(
-\frac{x^{2}}{2}\right)  }{1-\operatorname{erf}\left(  \frac{x}{\sqrt{2}%
}\right)  }\right]  ^{n+1}\left[  2\ln\left(  \frac{n+1}{1-\operatorname{erf}%
\left(  \frac{x}{\sqrt{2}}\right)  }\right)  \right]  ^{-\frac{1}{2}}.
\]

\begin{figure}[ptb]
\begin{center}
\rotatebox{270} {\resizebox{!}{5in}{\includegraphics{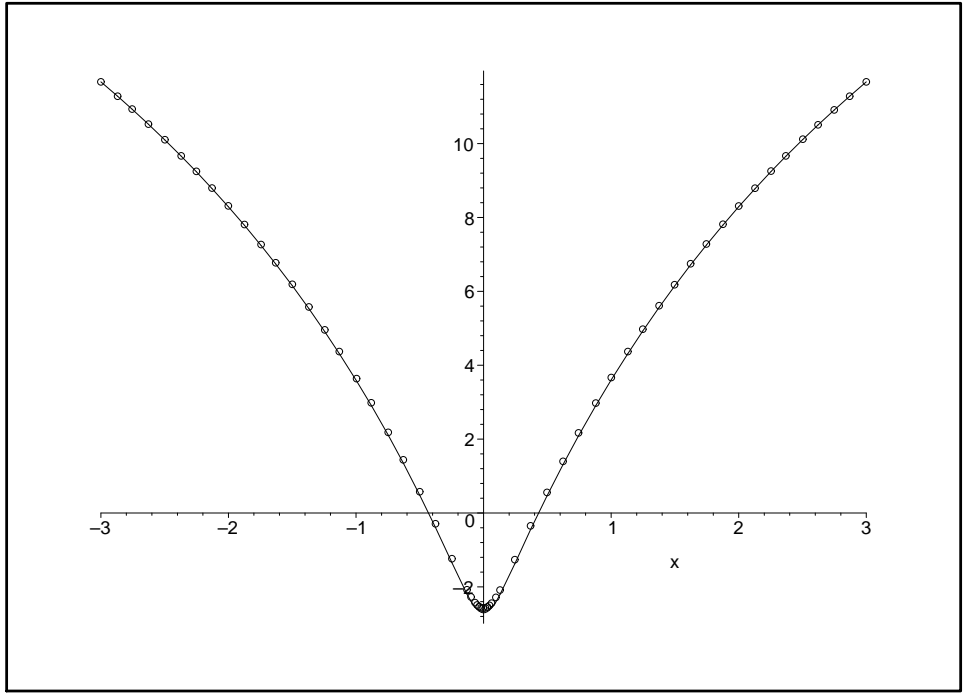}}}
\end{center}
\caption{A sketch of the exact (solid curve) and asymptotic (ooo) values of
$\ln\left[  \frac{P_{10}(x)}{10!}\right]  $.}%
\label{P10}%
\end{figure}

From (\ref{Pcnk}) we know that the polynomials $P_{n}(x)$ satisfy the
reflection formula
\begin{equation}
P_{n}(-x)=\left(  -1\right)  ^{n}P_{n}(x). \label{reflection}%
\end{equation}
Using (\ref{reflection}), we can extend (\ref{Pn5}) to the whole real line and
write%
\begin{equation}
P_{n}(x)\sim n!\left[  \Phi\left(  x,n\right)  +\left(  -1\right)  ^{n}%
\Phi\left(  -x,n\right)  \right]  ,\quad n\rightarrow\infty. \label{Pn6}%
\end{equation}
In Figure \ref{P10} we compare the values of $P_{10}(x)$ with the asymptotic
approximation (\ref{Pn6}).

We see that the approximation is very good, even for small values of $n.$ We summarize our results of this section in the following theorem.

\begin{theorem}
\label{theorem}Let the polynomials $P_{n}(x)$ be defined by%
\[
P_{n+1}(x)=P_{n}^{\prime}(x)+\left(  n+1\right)  xP_{n}(x),
\]
with $P_{0}(x)=1.$ Then, we have%
\begin{equation}
P_{n}(x)\sim P_{n}(x)\sim n!\left[  \Phi\left(  x,n\right)  +\left(
-1\right)  ^{n}\Phi\left(  -x,n\right)  \right]  ,\quad n\rightarrow\infty,
\label{Pasympt}%
\end{equation}
where
\begin{equation}
\Phi\left(  x,n\right)  =\left[  \sqrt{\frac{2}{\pi}}\frac{\exp\left(
-\frac{x^{2}}{2}\right)  }{1-\operatorname{erf}\left(  \frac{x}{\sqrt{2}%
}\right)  }\right]  ^{n+1}\left[  2\ln\left(  \frac{n+1}{1-\operatorname{erf}%
\left(  \frac{x}{\sqrt{2}}\right)  }\right)  \right]  ^{-\frac{1}{2}}.
\label{Phi}%
\end{equation}

\end{theorem}

\section{Higher derivatives of $\operatorname*{inverf}\left(  z\right)  $}

From (\ref{deriv2}) and (\ref{Pngn}), it follows that%
\begin{equation}
d_{n}=\frac{1}{\sqrt{2}}\left(  \sqrt{\frac{\pi}{2}}\right)  ^{n}%
P_{n-1}(0),\quad n=1,2,\ldots, \label{deriv3}%
\end{equation}
where $d_{n}$ was defined in (\ref{dn}). Using Theorem \ref{theorem} in
(\ref{deriv3}), we have%
\[
d_{n}\sim\frac{1}{\sqrt{2}}\left(  \sqrt{\frac{\pi}{2}}\right)  ^{n}%
\Phi\left(  0,n-1\right)  \left[  1+\left(  -1\right)  ^{n-1}\right]  ,
\]
as $n\rightarrow\infty.$ Using (\ref{Phi}), we obtain%
\begin{equation}
\frac{d_{n}}{n!}\sim\frac{1}{2n\sqrt{\ln(n)}}\left[  1+\left(  -1\right)
^{n-1}\right]  ,\quad n\rightarrow\infty. \label{deriv5}%
\end{equation}
Setting $n=2N+1$ in (\ref{deriv5}), we have%
\begin{equation}
\frac{d_{2N+1}}{\left(  2N+1\right)  !}\sim\frac{1}{\left(  2N+1\right)
\sqrt{\ln(2N+1)}},\quad N\rightarrow\infty. \label{derivodd}%
\end{equation}

\subsection{Numerical results}

\begin{figure}[ptb]
\begin{center}
\rotatebox{270} {\resizebox{!}{5in}{\includegraphics{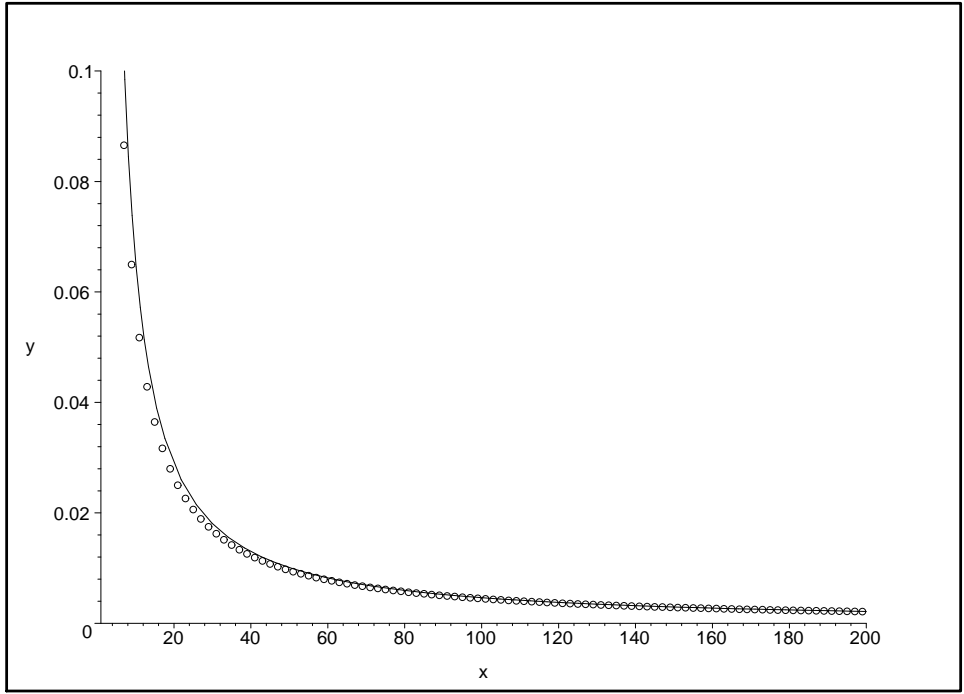}}}
\end{center}
\caption{A sketch of the exact (ooo) and asymptotic (solid curve) values of
$\frac{d_{2k+1}}{\left(  2k+1\right)  !}$.}%
\label{compare1}%
\end{figure}

In this section we demonstrate the accuracy of the approximation
(\ref{deriv5}) and construct a high order Taylor series for
$\operatorname*{inverf}\left(  x\right)  .$ In Figure \ref{compare1} we
compare the logarithm of the exact values of $\left.  \frac{d^{2n+1}%
}{dz^{2n+1}}\operatorname*{inverf}\left(  x\right)  \right\vert _{x=0}$ and
our asymptotic formula (\ref{deriv5}). We see that there is a very good
agreement, even for moderate values of $n$.

\begin{figure}[ptb]
\begin{center}
\rotatebox{270} {\resizebox{!}{5in}{\includegraphics{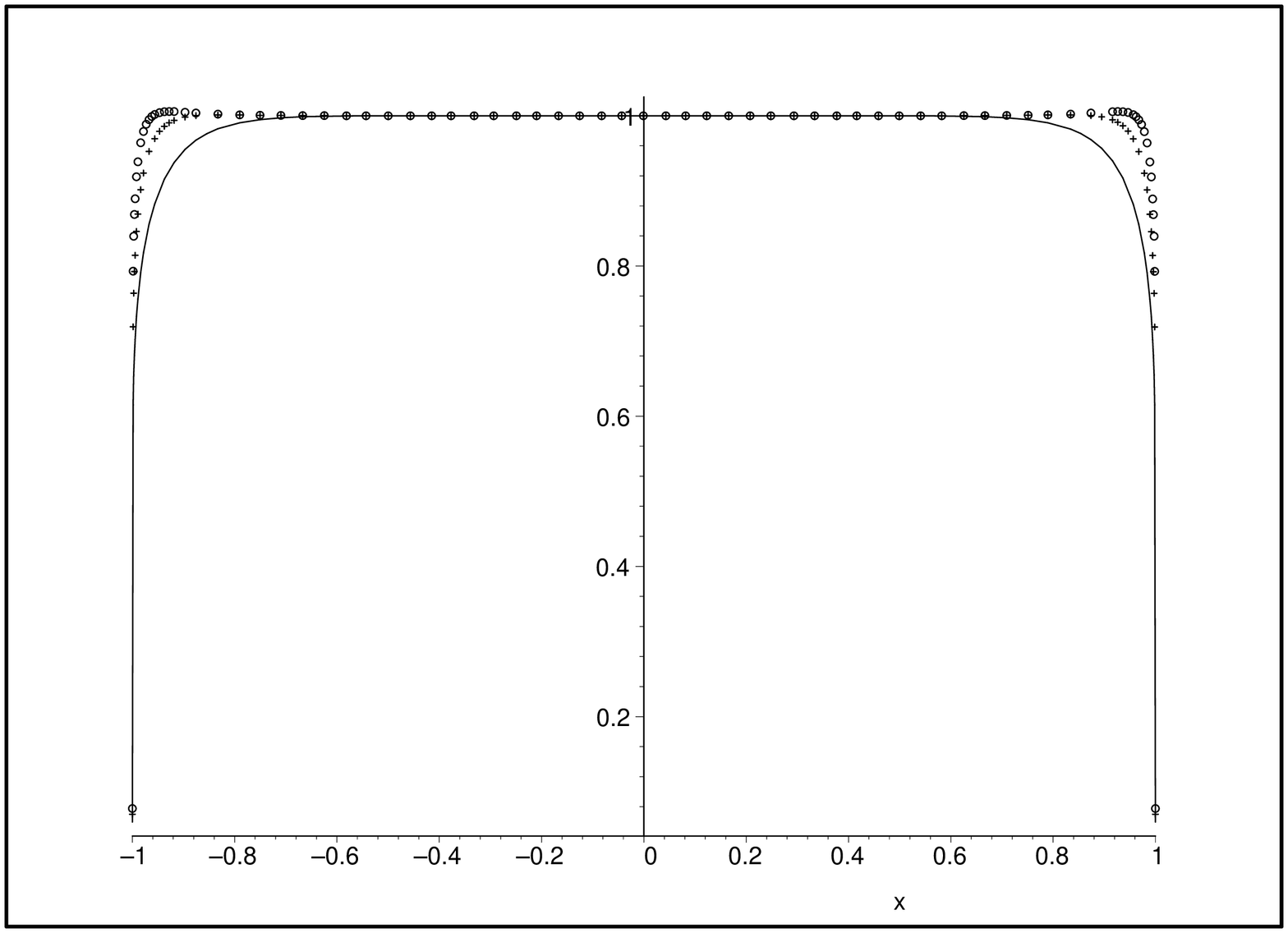}}}
\end{center}
\caption{A sketch of $\frac{T_{9}(x)}{\operatorname*{inverf}\left(  x\right)
}$ (solid curve), $\frac{T_{9}(x)+R_{10}(x)}{\operatorname*{inverf}\left(
x\right)  }$ (+++) and $\frac{T_{9}(x)+R_{20}(x)}{\operatorname*{inverf}%
\left(  x\right)  }$ (ooo).}%
\label{compare2}%
\end{figure}

Using (\ref{recurrence}), we compute the exact values
\[
d_{1}=\frac{1}{2}\pi^{\frac{1}{2}},\quad d_{3}=\frac{1}{4}\pi^{\frac{3}{2}%
},\quad d_{5}=\frac{7}{8}\pi^{\frac{5}{2}},\quad d_{7}=\frac{127}{16}%
\pi^{\frac{7}{2}},\quad d_{9}=\frac{4369}{32}\pi^{\frac{9}{2}}%
\]
and form the polynomial Taylor approximation%
\[
T_{9}(x)=\sum\limits_{k=0}^{4}d_{2k+1}\frac{x^{2k+1}}{\left(  2k+1\right)
!}.
\]
In Figure \ref{compare2} we graph $\frac{T_{9}(x)}{\operatorname*{inverf}%
\left(  x\right)  }$ \ and \ $\frac{T_{9}(x)+R_{N}(x)}{\operatorname*{inverf}%
\left(  x\right)  },$ for $N=10,20,$ where
\begin{equation}
R_{N}(x)=\sum\limits_{k=5}^{N}\frac{x^{2k+1}}{\left(  2N+1\right)  \sqrt
{\ln(2N+1)}},\quad N=5,6,\ldots. \label{T1}%
\end{equation}

\begin{figure}[ptb]
\begin{center}
\rotatebox{270} {\resizebox{!}{5in}{\includegraphics{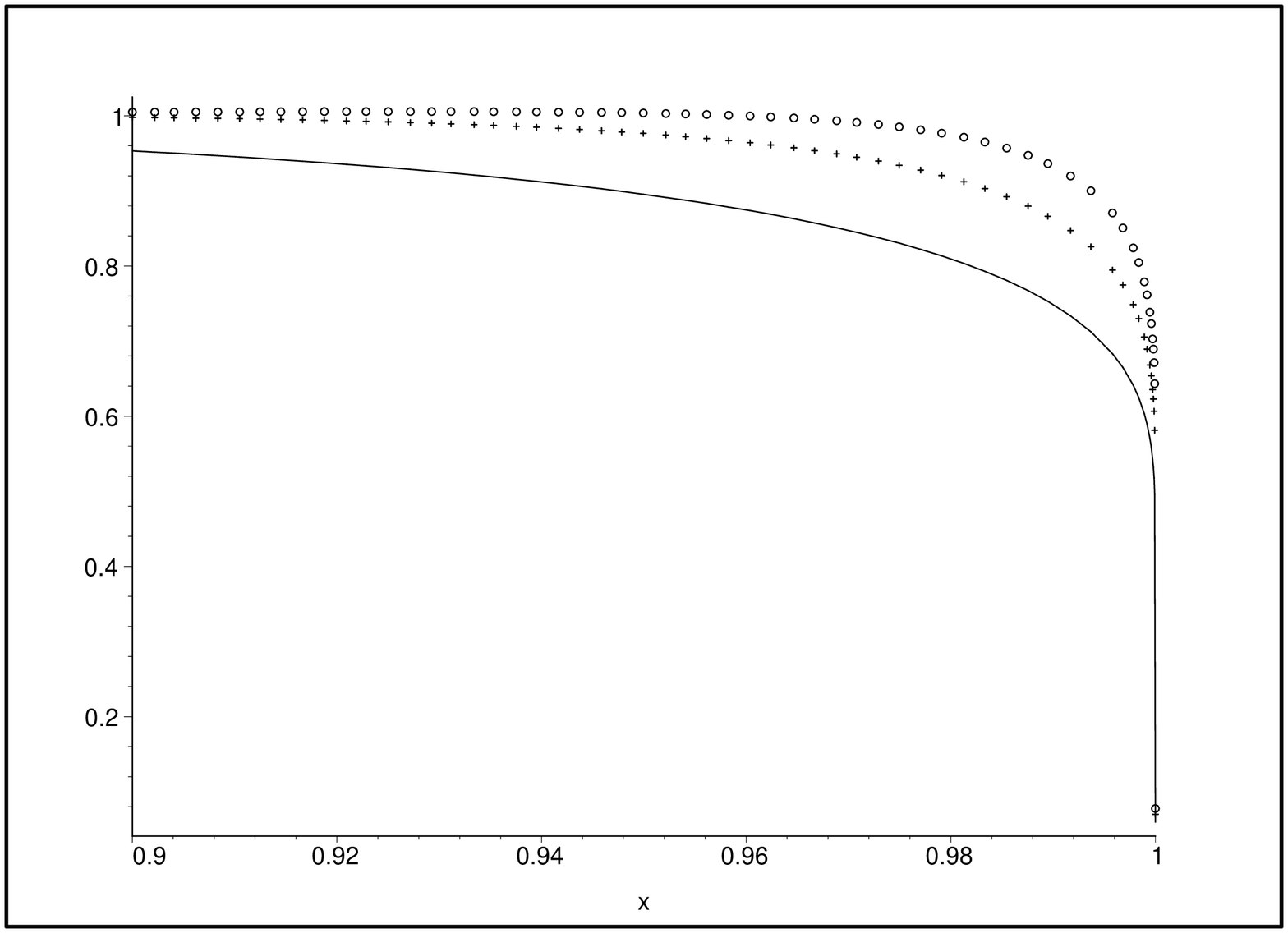}}}
\end{center}
\caption{A sketch of $\frac{T_{9}(x)}{\operatorname*{inverf}\left(  x\right)
}$ (solid curve), $\frac{T_{9}(x)+R_{10}(x)}{\operatorname*{inverf}\left(
x\right)  }$ (+++) and $\frac{T_{9}(x)+R_{20}(x)}{\operatorname*{inverf}%
\left(  x\right)  }$ (ooo).}%
\label{compare3}%
\end{figure}

The functions are virtually identical in most of the interval $\left(
-1,1\right)  $ except for values close to $x=\pm1.$ We show the differences in
detail in Figure \ref{compare3}. Clearly, the additional terms in $R_{20}(x)$
give a far better approximation for $x\simeq1.$

In the table below we compute the exact value of and optimal asymptotic
approximation to $\operatorname*{inverf}\left(  x\right)  $ for some $x$:%
\[%
\begin{tabular}
[c]{|c|c|c|c|}\hline
$x$ & $\operatorname*{inverf}\left(  x\right)  $ & $T_{9}(x)+R_{N}(x)$ &
$N$\\\hline
$0.7$ & $.732869$ & $.732751$ & $6$\\\hline
$0.8$ & $.906194$ & $.905545$ & $7$\\\hline
$0.9$ & $1.16309$ & $1.16274$ & $11$\\\hline
$0.99$ & $1.82139$ & $1.82121$ & $57$\\\hline
$0.999$ & $2.32675$ & $2.32676$ & $423$\\\hline
$0.9999$ & $2.75106$ & $2.75105$ & $3685$\\\hline
\end{tabular}
\ \ .
\]
Clearly, (\ref{T1}) is still valid for $x\rightarrow1,$ but at the cost of
having to compute many terms in the sum. In this region it is better to use
the formula \cite{MR1986919}%
\[
\operatorname*{inverf}\left(  x\right)  \sim\sqrt{\frac{1}{2}%
\operatorname*{LW}\left[  \frac{2}{\pi\left(  x-1\right)  ^{2}}\right]
},\quad x\rightarrow1^{-},
\]
where $\operatorname*{LW}(\cdot)$ denotes the Lambert-W function
\cite{MR1414285}, which satisfies%
\[
\operatorname*{LW}(x)\exp\left[  \operatorname*{LW}(x)\right]  =x.
\]

\bibliographystyle{amsalpha}
\newcommand{\etalchar}[1]{$^{#1}$}
\providecommand{\bysame}{\leavevmode\hbox to3em{\hrulefill}\thinspace}
\providecommand{\MR}{\relax\ifhmode\unskip\space\fi MR }
\providecommand{\MRhref}[2]{%
  \href{http://www.ams.org/mathscinet-getitem?mr=#1}{#2}
}
\providecommand{\href}[2]{#2}

\end{document}